\documentclass[letterpaper,11pt]{article}
\usepackage{amssymb,amsmath,amsthm,turnstile,txfonts}
\usepackage{color}
\usepackage{graphicx}
\usepackage{epstopdf}
\usepackage{url}
\usepackage{caption}
\usepackage{multirow}
\usepackage[left=1in,top=1in,right=1in,bottom=1in,letterpaper]{geometry}
\usepackage{setspace}
\usepackage[algo2e,linesnumbered, vlined,ruled]{algorithm2e}
\usepackage{wrapfig}
\usepackage{float}

\newcommand\myfigure[1]{%
\medskip\noindent\begin{minipage}{\columnwidth}
\centering%
\end{minipage}\medskip}

\def\g{\textbf{g}}

\def\w{\textbf{w}}
\def\x{\textbf{x}}

\def\M{\mathcal{M}}

\def\E{\mathcal{E}}

\def\W{\mathcal{W}}

\def\S{\mathcal{S}}

\def\RR{\mathbb{R}}

\def\Fun{\mathrm{C}^{\infty}}

\DeclareMathOperator{\tr}{tr}

\DeclareMathOperator{\argmin}{arg min}

\newtheorem{theorem}{{\bf Theorem}}
\newtheorem{proposition}{{\bf Proposition}}

\makeatletter
\@addtoreset{proofpart}{theorem}
\makeatother

\title{Nonisometric Surface Registration via \\
Conformal Laplace-Beltrami Basis Pursuit
}
\author{ 
Stefan C. Schonsheck
    \thanks{Department of Mathematics, Rensselaer Polytechnic Institute, Troy, NY 12180, U.S.A. ({\tt schon@rpi.edu}). }
\and Michael M. Bronstein
     \thanks{Institute of Computational Science; Universit� della Svizzera Italiana,  ({\tt michael.bronstein@gmail.com}).}
\and Rongjie Lai
\thanks{Department of Mathematics, Rensselaer Polytechnic Institute, Troy, NY 12180,
         U.S.A. ({\tt lair@rpi.edu}). }
        }
\date{}

\begin{document}
\maketitle

\begin{abstract} Surface registration is one of the most fundamental problems in geometry processing. Many approaches have been developed to tackle this problem in cases where the surfaces are nearly isometric. However, it is much more challenging  to compute correspondence between surfaces which are intrinsically less similar. In this paper, we propose a variational model to align the Laplace-Beltrami (LB) eigensytems of two non-isometric genus zero shapes via conformal deformations. This method enables us compute to geometric meaningful point-to-point maps between non-isometric shapes. Our model is based on a novel basis pursuit scheme whereby we simultaneously compute a conformal deformation of a 'target shape' and its deformed LB eigensytem. We solve the model using an proximal alternating minimization algorithm hybridized with the augmented Lagrangian method which produces accurate correspondences given only a few landmark points. We also propose a reinitialization scheme to overcome some of the difficulties caused by the non-convexity of the variational problem. Intensive numerical experiments illustrate the effectiveness and robustness of the proposed method to handle non-isometric surfaces with large deformation with respect to both noise on the underlying manifolds and errors within the given landmarks. 
\end{abstract}

\section{Introduction}\label{sec:Intro}

The computation of meaningful point-to-point mappings between pairs of manifolds lies at the heart of many shape analysis tasks. It is crucial to have robust methods to compute dense correspondences between two or more shapes in different applications including shape matching, label transfer, animation and recognition~\cite{siddiqi1998area,kraevoy2004cross,reuter2006laplace,heimann2009statistical,van2011survey,ovsjanikov2012functional}. In cases where shapes are very similar (isometric or nearly isometric), there are many approaches for computing such correspondences \cite{elad2003bending,gu2004genus,bronstein2006efficient,aubry2011wave,kim2011blended,kovnatsky2013coupled,lai2017multi,ovsjanikov2012functional,shtern2015spectral,shtern2016fast}. However, it is still challenging to compute accurate correspondences when the shapes are far away from near isometry.

One of the key challenges in largely deformed non-isometric shape matching is that the intrinsic features of the two shapes are not similar enough for standard techniques to recognize their similarity. For example, when computing the correspondence between human faces, it is not particularly difficult to geometrically characterize the structure of a `nose'. However, similar techniques can not work well to compute a map between a horse and an elephant face since these two surfaces have many largely deformed local structures including the drastic difference between the trunk of the elephant and the nose of the horse. Because of this, it is crucial to develop new methods to adaptively characterize large deformations on surfaces. 

The LB eigensystem is a ubiquitous tool for 3D shape analysis (see \cite{berard1994embedding,reuter2006laplace,Levy06,Rustamov:2007,sun2009concise,Bronstein:2010CVPR,Lai:2010CVPR,Raviv2011CVPR,aubry2011wave,lai2011automated,shi2013cortical,shi2014metric,shtern2015spectral} and references therein). It is invariant under isometric transformations and intrinsically characterizes the local and global geometry of manifolds through its eigensystem up to an isometry. In principle, the LB eigensystem reduces infinite dimensional nonlinear isomorphism ambiguities between two isometric shapes to a linear transformation group between two LB eigensystems. This linear transform is necessary due to the possible sign or sub-eigenspace ambiguity of LB eigensystems~\cite{lai2017multi}. Additionally, similar shapes often have similar eigensystems which allows for joint analysis of similar shapes their spectral properties~\cite{ovsjanikov2012functional}. However, when the deformation between two shapes is far from an isometry, the large dissimilarity between LB eigensystems of two shapes is the major bottleneck to adapt existing spectral geometry approach to conduct registration. 

A natural idea to extend spectral geometry methods to register non-isometric surfaces is to deform the metric of the target surface to the source one such that two surfaces share similar LB eigensystem after deformation. However,  the amount of imposed deformation is essentially relied on the correspondence between these two surfaces. In this work, we proposed to simultaneously compute the deformation and learn features for registration. Mathematically, one way to characterize these deformations is through the conformal factor. It is well known that there exists a conformal mapping between any two genus zero surfaces~\cite{jost2008riemannian}. Rather than reconstruct conformally deformed surfaces and related conformal map, we exploit a fundamental link between the conformal factor and the LB eigesystem by considering the conformally deformed LB eigensystem. This allows us to compute a new basis on the target surface to align the naturally defined LB eigensystem on the source surface. It leads to a variational method for non-isometric shape matching which enables us to overcome the natural ambiguities of the LB eigensystem, to align the bases of non-isometric shapes and to avoid the direct computation of conformal maps. 

Numerically, we solve our model using a proximal alternating minimization (PAM) method~\cite{attouch2010proximal} hybridized with the augmented Lagrangian method~\cite{glowinski1989augmented}. The method is iteratively composed of a curvilinear search method on orthogonality constrained manifold~\cite{wen2013feasible} in one direction to compute the conformally deformed LB eigenfunctions and  the BFGS~\cite{bazaraa2013nonlinear} method for the other direction to compute the conformal factor. Theoretically, we can guarantee the local convergence of the proposed algorithm since the objective function and constraints satisfy the necessary Kurdyka-Lojasiewicz (KL) condition~\cite{attouch2010proximal}. Numerical results on largely deformed test problems, including horse-to-elephant and Faust benchmark database~\cite{bogo2014faust}, validate the effectiveness and robustness of our method.

\paragraph{Related Works.}
 A large number of 3D nonrigid shape matching approaches are based on analysis of the LB eigen system~(see \cite{reuter2005laplace,reuter2006laplace,Levy06,Rustamov:2007,Bronstein:2010CVPR,ovsjanikov2012functional,kovnatsky2013coupled,Raviv:2014Affine,shtern2015spectral,lai2017multi} and reference therein). The LB eigensytem is intrinsic and invariant to isomorphism, and also characterizes the local and global geometry of a manifold. This makes it ideal for many shape processing tasks and many early works in the field involve directly comparing the LB spectrum of the shapes to determine how alike shapes are \cite{reuter2005laplace,reuter2006laplace,Levy06}. More recently, the general concept of functional maps \cite{ovsjanikov2012functional} has played a central role in many new methods which have allowed for the formulation of accurate correspondence maps. This technique essentially reduces the non-linear transform between two shapes to a linear transform between their eigensystems. In general, these techniques work for well for isometric and near isometric cases, but can not produce satisfactory results when the LB eigensystems of shapes are very dissimilar. This occurs when the deformation between shapes is far from an isometry. To overcome this, the concept of coupled bases (also known as joint-diagonalization) was introduced for shape processing tasks in \cite{kovnatsky2013coupled}. In this work the authors propose a variational model to define a shared basis for a pair of shapes which is `nearly harmonic' on one shape and 'similar' to the natural LB basis on the other. This joint optimization allows for much more accurate correspondence maps, but does not characterize the underlying deformations which lie at the heart of the non-isometric shape matching problem.

Conformal maps have been widely applied to various shape processing tasks in order to characterize these deformations \cite{Hurdal:NeuroImage2000,Haker:TVCG2000,gu2004genus,Springborn:SIGGRAPH08}.  
In one of the first works to combine spectral and deformation based approaches, \cite{shi2011conformal} presents a scheme to find optimal conformal deformation to align two shapes in the embedded LB Space. Additionally, the authors present a general framework for computing LB eigensytems of conformaly deformed surfaces as well as several other imported related quantities. Continuing on this line of work in \cite{kao2017maximization}, the authors use the LB eigenvalues as a tool to guide conformal deformations. Using derivatives of the the LB eigenvalues, they compute optimal conformal metrics which approximate conformal and topological eigenvalues. In our work, we use the spectral coefficients of known features to guide the deformation, so rather than align the eigenvalues we align the eigenfunctions. This allows us to avoid the subspace ambiguity of the the LB eigensystem and computational errors in calculating high frequency eigenvalues. 

\paragraph{Major Contributions.}
We introduce a novel variational basis pursuit model for computing non-isometric shape correspondences via conformal deformation of the LB eigensystem. This model enhances spectral approaches from handling nearly isometric surface registration to tackling surfaces with large deformed metrics. It naturally combines the conformal deformation to the LB eigensystem and simultaneously computes surface deformations and LB eigenbasis which also automatically overcomes the ambiguities of LB eigensystems in surface registration. 
We also propose a numerical scheme to solve the variational model with local convergence guarantee. Additionally, we introduce a reinitializaiton scheme to help tackle local minima and improve the quality of the computed bases. This algorithm successfully handles non-isomorphic shape correspondence problems given only a few landmarks and is shown to be robust to noise and perturbations of landmarks.

The rest of this paper is organized as follows: In section \ref{sec:Backgroud}, we review the theoretical background of conformal deformation of LB eigensystem and functional maps. After that, we propose the variational basis pursuit model for conformal deformation of the LB eigensystem in section \ref{sec:LBBasisPursuit}. In section \ref{sec:Algs}, we discretize the model and develop an optimization scheme based on PAM to solve the variational problem. Section \ref{sec:Discussion} is further devoted to discuss a few details of the model and a reinitialization scheme to improve our numerical solver. In section \ref{sec:Results}, numerical results on several data sets are presented to show that the model accurately produces point-to-point mappings on non-isometric manifolds with large deformation given only a few landmark points. We also show that our approach is robust to both noise in the underlying manifolds and inaccuracies in the initial landmarks. Furthermore, we test the model to a benchmark data based to show its effectiveness. Lastly, we conclude the paper in Section~\ref{sec:conclusion}.

\section{Mathematical Background}\label{sec:Backgroud}
In this section, we discuss the mathematical background of the proposed method. We first review a few key properties of the LB eigensystem of a Riemannian surface and discuss its conformal deformations with respect to deformations of the Riemannian surface metric \cite{chavel1984eigenvalues,jost2008riemannian}. After this, we review the functional maps framework in \cite{ovsjanikov2012functional} which will be closely related to our work.
\subsection{Conformal deformation of LB eigensystem on Riemannian Surfaces}
\label{subsec:LB}

Given a closed Riemannian surface $(\M,g)$, its LB operator in a given local coordinate system, $\{x_i\}_{i=1,2}$,  is defined as \cite{chavel1984eigenvalues,jost2008riemannian}:
\begin{equation}
\Delta_{g}\phi=\frac{1}{\sqrt{G}}\sum_{i=1}^{2}\frac{\partial}{\partial
x_i}(\sqrt{G}\sum_{j=1}^{2}g^{ij}\frac{\partial \phi}{\partial x_j})
\label{def:LB}
\end{equation}
where $(g^{ij})$ is the inverse of the metric matrix $g = (g_{ij})$

and $G=\det(g_{ij})$. The LB operator is self-adjoint and elliptic, therefore it has a discrete spectrum. We denote the eigenvalues of $-\Delta_{g}$ as $0=\lambda_0 < \lambda_1 \leq \lambda_2 \leq\cdots$ with the corresponding eigenfunctions $\phi_0, \phi_1,\phi_2,\cdots$ satisfying:
\begin{equation} 
-\Delta_g(x) \phi_i(x) =\lambda_i\phi_i(x), \quad \text{and} \quad 
\int_{\M}\phi_i(x) \phi_j(x) ~\mathrm{d}vol_{g}(x)= \delta_{ij}, \quad i,j = 0,1,2,\cdots 
\end{equation}
where $\mathrm{d}vol_{g}(x)$ is the area element on $\M$ with respect to $g$. It is well-known that $\Phi = \{\phi_n~|~ n=0, 1, 2, \cdots\}$ forms an orthonormal basis for the real-valued, smooth function space $\Fun (\M,\RR)$ on the manifold $(\M,g)$. This basis can be viewed as a generalization of the Fourier basis from flat space to a differentiable manifold. 
The LB eigensystem is invariant under both rigid and nonrigid isometric transformations,
and it uniquely determines a manifold up to isometry \cite{berard1994embedding}.

In differential geometry, a conformal map is one which preserves angles locally. Formally, a conformal map preserves the first fundamental form up to a positive scaling factor. Given two manifolds $(\M_1,g_1)$ and $(\M_2,g_2)$, a map $F:(\M_1,g_1) \rightarrow (\M_2,g_2)$ is conformal if and only if the pullback $F^*(g_2) = w^2 g_1$ with a positive function $w^2$ (written this way to emphasize positivity). A {\it conformal deformation} of a surface is a transformation which changes the local metric by a positive scaling factor. A well known result in conformal geometry is that there exists a conformal map between any two genus zero surfaces \cite{jost2008riemannian}. 

Given a closed surface $(\M,g)$ with conformal deformation $w^2$, the LB eigensystem of the deformed manifold $(\M, w^2 g)$ can be viewed as a weighted LB eigensystem on the original surface $(\M,g)$. This simple fact intrinsically links the LB eigensystem of the deformed manifold to a weighed LB eigensystem on the original manifold. It allows us to compute the LB eigensytem of the conformally deformed manifold without explicitly reconstructing its embedding or coordinates. This also relates information about the local deformation and global eigensytem and later becomes the cornerstone of our approach. Formally, we have: 
\begin{proposition}
\label{prop:LBconformal}
Let $\{\phi_n^{w^2},\lambda^{w^2}_n\}_{n=1}^\infty$ be a LB eigensystem of a conformally deformed surface $(\M,w^2 g)$, then $\{\phi_n^{w^2},\lambda^{w^2}_n\}_{n=1}^\infty$ is equivalent to the following weighted LB eigensystem on $(\M,g)$:
\begin{eqnarray}\label{eqn:ConformalLBeigs}
-\Delta_{g}\phi_i(x) = \lambda  w^2(x)\phi_i(x), \qquad
\int_{\M}\phi_i(x)\phi_j(x) w^2(x)~ c(x)= \delta_{ij}
\end{eqnarray}
\end{proposition}
\begin{proof} This is because  
\begin{align*}
\Delta_{w^2g} \phi
= \frac{1}{w^2\sqrt{G}}\sum_{i=1}^{2}\frac{\partial}{\partial
x_i}(w^2\sqrt{G}\sum_{j=1}^{2}w^{-2}g^{ij}\frac{\partial \phi}{\partial x_j}) = w^{-2}\Delta_{g} \phi
\end{align*}
Hence the eigen problem: $-\Delta_{w^2g} \phi = \lambda \phi$ is equivalent to $-\Delta_{g}\phi = \lambda  w^2\phi$. In addition, $\mathrm{d}vol_{w^2 g} = w^2\ \mathrm{d}vol_{g}$
\end{proof}

The problem of finding the LB eigensystem of a Riemannian manifold is equivalent to finding an orthonormal set of functions $\Phi = \{\phi_i\}$ which have minimal harmonic energy on the surface. From the above proposition, the LB eigensystem of a conformally deformed manifold $(\M,w^2 g)$ can be formulated as the following variational problem: 
\begin{equation}\label{defeig}
\arg\min_{\Phi=\{\phi_i\}} \sum_i  \int_{\M} ||\nabla_{\M} \phi_i(x) ||^2~\mathrm{d}vol_{g}(x),  \quad \text{s.t.}\quad \int_M \phi_i(x) \phi_j(x) w^2(x) ~\mathrm{d}vol_{g}(x) = \delta_{ij}
\end{equation}

\subsection{Functional Maps}
\label{funmaps}
Functional maps were introduced in~\cite{ovsjanikov2012functional} for isometric and nearly isometric shape correspondence. This method has been shown a very effective tool for various shape processing tasks \cite{ovsjanikov2012functional,kovnatsky2013coupled,rodola2015point}. Here we provide a basic overview of their framework.
Consider Riemannian surfaces $(\M_1,g_1)$ and $(\M_2,g_2)$, a smooth bijection $F:\M_1 \rightarrow \M_2$ induces a linear transformation between functional spaces of these two manifolds as:
$$
F_T: \Fun(\M_1,\RR) \rightarrow  \Fun(\M_2,\RR),\quad 
f \mapsto f\circ F^{-1}$$
Instead of computing surface map $F$, the crucial idea of functional map is to  compute the linear map $F_T$ between these two functional spaces. After that, the desired surface map can be encoded by considering images of indicator functions under $F_T$.

Finding a functional map, $F_T$, associated with a map $F$ is equivalent to finding the matrix representation of $F_T$ under a fixed orthonormal basis $\{\phi_i\}$ of $\Fun(\M_1,\RR)$ and a fixed orthonormal basis $\{\psi_i\}$ of $\Fun(\M_2,\RR)$, respectively.  Namely, if we write $F_T(\phi_i) = \sum_{j}c_{ji} \psi_j$, then any two given corresponding functions $ f =  \sum_i f_i \phi_i$ and $g = \sum_j g_j \psi_j$ under $F_T$ can be represented using $C = (c_{ij})$ as: 
$$ F_T(f) = g \Leftrightarrow F_T\Big(\sum_i f_i \phi_i\Big) 
=\sum_i f_i F_T ( \phi_i) =\sum_i f_i \sum_j c_{ji} \psi_j = \sum_j g_j \psi_j 
\Leftrightarrow \sum_i c_{ji} f_i = g_j.$$
Each entry of the matrix $c_{ij}$ can be found by finding the $j^{th}$ coefficient of $F_T(\phi_i)$ expressed in the $\{\psi_i\}$ coordinate system, i.e. $c_{ji} = \langle F_T(\phi_i), \psi_j \rangle_{g_2}$. In practice, one can use two finite sets of orthonormal functions to approximate $\Fun(\M_1,\RR)$ and $\Fun(\M_2,\RR)$, thus the functional map can be approximated by a finite dimensional matrix. For instance, the first $N$ eigenfunctions of the LB eigensystem is one common choice of such a basis. Then, the problem of finding the transformation $F_T$ can be approximated by the problem of seeking a finite dimension matrix $C$. As long as $C$ is computed, the desired map $F$ can be computed through $C$ operating on indicator functions. 

\section{Conformal LB Basis Pursuit for Nonisometric Surface Registration}
\label{sec:LBBasisPursuit}
In this section, we propose a LB basis pursuit model for non-isometric surface registration. On the target surface $\M_2$, the model simultaneously finds a conformal deformation and a conformally deformed LB eigensystem so that the coefficients of the corresponding feature functions expressed on the deformed LB eigensystem of $\M_2$ are the same as the coefficients on the fixed source surface $\M_1$. 

\subsection{Variational PDE model}

Given two non-isometric genus-0 closed Riemannian surfaces $(\M_1,g_1), (\M_2,g_2)$, we aim at finding a geometrically meaningful correspondence between these two surfaces. In the case that $\M_1$ and $\M_2$ are nearly isometric, there are many successful methods to constructing maps between $\M_1$ and $\M_2$ by comparing their isometric invariant features. Using spectral descriptors from solutions of the LB eigensystem on manifolds is a common way of constructing such descriptors~\cite{reuter2006laplace,Levy06,Vallet:2008CGF,Shi:08a,Bronstein:2010CVPR}. As extensions, some other descriptors such as Heat kernel signature~\cite{sun2009concise}, wave kernel signature~\cite{aubry2011wave} and optimal spectral descriptors~\cite{litman2014learning} have also been proposed in the literature. However, most of the existing methods consider the construction of descriptors for nearly isometric manifolds. Registration methods based on the existing LB spectral descriptors can not provide satisfactory results for constructing correspondence between two non-isometric surfaces as their eigensystems are possibly quite far apart. 

We propose to overcome the limitation of the LB spectral descriptors for largely deformed non-isometric shape registration by considering a continuous deformation of the LB spectral descriptors. Intuitively, given two non-isometric shapes $(\M_1,g_1)$ and $(\M_2,g_2)$, our idea is to deform the metric of $(\M_2,g_2)$ such that the deformed surface is isometrically the same as $(\M_1,g_1)$. Then the LB spectral descriptors can be applied as in isometric shape matching. However, it is challenging to find an appropriate deformation as the accurate amount of deformation on each local region of $\M_2$ depends exactly on an accurate correspondence which is precisely the problem we would like to solve.

To handle this challenge, we propose to simultaneously find an optimal correspondence and an optimal deformation. More specifically, by fixing the LB eigensystem $\{\Phi, \Lambda\}$ of $(\M_1, g_1)$, we seek a map $T:\M_1 \rightarrow \M_2$ and a conformal factor $w^2:\M_2 \rightarrow \RR^{+}$ such that the LB eigensystem $\{\Phi, \Lambda\}$ of $(\M_1,g_1)$ can be aligned to the LB eigensystem $\{\Psi, \Theta\}$ of $(\M_2, w^2 g_2)$ via $T$. 
This problem can be written as the following variational PDE problem:
\begin{equation}\begin{split}
\label{eqn:map1}
(T^*,w^*, \Psi^*) &=
\text{argmin}_{T, w,\Psi = \{\psi_i\}_{i=1}^N } \sum_{i=1}^{N} \int_{\M_1} \|\phi_i - \psi_i\circ T\|^2 ~\mathrm{d}M_1  + \frac{1}{2} \sum_{i=1}^N \int_{\M_2} \|\nabla_{\M_2} \psi_i \|^2 ~\mathrm{d}\M_2,  \\
 &\hspace{4cm} \text{s.t.}\quad  \int_{\M_2} \psi_i \psi_j \ w^2 ~\mathrm{d}\M_2  = \delta_{ij}
\end{split}
\end{equation}
where $\mathrm{d}\M_1 =dvol_{g_1}, \mathrm{d}\M_2 = dvol_{g_2} $ and $w^2 \mathrm{d} \M_2 = dvol_{w^2 g_2}$.The first term measures the alignment of two bases as the correct correspondence should map one LB eigensystem to another one, and the second term solves the first $N$ LB eigenfunctions $\{\psi_i\}$ for the deformed manifold $(\M_2,w^2 g_2)$ due to the variational problem~\eqref{defeig}. 
Existence of a solution to this variational problem \eqref{eqn:map1} is guaranteed as any two genus-0 surfaces are conformally equivalent and the LB operator is invariant under isometric transformations.

Computationally, the numerical search for $T$ in the mapping space is usually very time consuming. Inspired by the idea of functional maps \cite{ovsjanikov2012functional} and the coupled quasi-harmonic bases \cite{kovnatsky2013coupled}, we choose to represent $T$ in the functional space. Instead of finding $T$ directly, we look for a basis $\Psi = \psi_i\circ T = F_T(\psi_i)$ which is nearly harmonic on $(\M_2,w^2 g_\M)$ and represents the corresponding features with the same coefficients as $\Phi$ does.  
More precisely, given a set of corresponding features $F=\{f_1,\cdots, f_k\} $ on $\M_1$ and $G = \{g_1, \cdots, g_k\}$ on $\M_2$, such that $f_i(x) = g_i(y)$ if $x$ and $y$ are corresponding points on $\M_1$ and $\M_2$, we can replace the direct measurement of the basis alignment term with a coefficient matching term. 
That is, instead of measuring the alignment of $\Psi$ and $\Phi$ via $T$, we measure how closely the coefficients for $G$ in the computed basis $\Psi$ match the coefficients for $F$ in the fixed LB basis $\Phi$. Formally, we measure the coefficient alignment by constructing a matrix of the coefficients in for $F$ in $\Phi$ and for $G$ in $\Psi$ so that the $ij^{th}$ term represents the coefficient for the $i^{th}$ corresponding function in the $j^{th}$ basis and computing their difference under the Frobenius norm. With this in mind, we propose the following model:
\begin{equation}
\label{eqn:map2}\begin{split}
(w^*,\Psi^*) &= \arg\min_{ w,\Psi} 
\frac{r_1}{2}\| \langle F, \Phi \rangle_{g_1} -  \langle G, \Psi \rangle_{w^2 g_2} \|_F^2 + \frac{r_2}{2} \sum_{i=1}^N \int_{\M_2} \|\nabla_{\M_2} \psi_i \|^2 \mathrm{d}\M_2, \\
&\hspace{3cm} \text{s.t.}\quad \int_{\M_2} \psi_i \psi_j \ w^2 \mathrm{d}\M_2  = \delta_{ij}
\end{split}
\end{equation}
where: 
$$\langle F, \Phi \rangle_{g_1} = \Big(\int_{\M_1} f_i \phi_j \ \mathrm{d}\M_1 \Big)_{i,j = 1,2,\dots,k} \quad 
\text{ and }\quad  \langle G, \Psi \rangle_{w^2 g_2} = \Big(\int_{\M_2} g_i \psi_j\  w^2 \mathrm{d}\M_2 \Big)_{i,j = 1,2,\dots,k}.$$ 
In practice we use indicator functions for $F$ and $G$, but heat signatures~\cite{sun2009concise}, wave kernel signatures~\cite{aubry2011wave}, or any other corresponding functions will also work.  Once $\Psi^* = \{\psi^*_1, \cdots, \psi^*_{M_2} \}$ is obtained, we can easily compute the functional map as 
\begin{equation}
\label{eqn:funmap}
F_T : C^\infty(\M_1) \rightarrow C^\infty(\M_2), \qquad F_T(h) = \sum_{i=1}  \Big(\int_{\M_1} h \phi_i \ \mathrm{d}vol_{g_1}\Big)~  \psi^T_i.
\end{equation} 

The main advantage of this model over previous existing methods for shape correspondence is that we are able to employ much more of the information encoded in the differential structures of $\M_1$ and $\M_2$ in our algorithm by combining the spectral descriptors and local deformations. This additional flexibility enables us to compute correspondences between largely deformed shapes. Information about the conformal deformation of the metric allows us to find a harmonic basis on the deformed shape, meanwhile information about the alignment of the functional spaces guides our calculation of the conformal deformation. Furthermore the additional constraint of the feature alignment overcomes ambiguity casued by the fact that there is no unique conformal deformation between any two genus zero surfaces.
To the best of our knowledge, the link between the conformal factor and deformed LB basis has not been exploited in such a way. Previous works have used only the conformal factor \cite{gu2004genus,kim2011blended} or only the functional space \cite{ovsjanikov2012functional,kovnatsky2013coupled} as stand alone tools rather than in concert as we present here. 

\subsection{Regularization and Area Constraint}

We add harmonic energy term to smooth the conformal deformation and regularize the problem. This can both increase the speed of the algorithm and improve the quality of the map, both in terms of the geodesic errors of the final correspondence, and the accuracy of the resulting conformal factor. This is particularly helpful to handle  deformations between the shapes which are far from isometry and to reduce the required number of features. Rather than smooth the conformal factor $w^2$ directly, we instead add the harmonic energy of $w$ to our objective function. Using $w$ instead of $w^2$ allows for easier analytic computation of the derivatives and a more efficient algorithm. In cases where the deformations are likely to be highly localized, this term may be omitted.

Lastly, we add an area preservation constraint to our model. That is, we would like the final deformed shape to be of the same size as the one we are matching it to. To enforce this, we mandate that the deformed manifold have the same surface area as the original manifold. This eliminates any scaling ambiguity. Then the final version of our model can be stated as:
\begin{equation} \label{eqn:finalmodel}
\begin{split}
(w^*,\Psi^*) =& \arg\min_{ w,\Psi=\{\psi_i\}_{i=1}^N} 
\frac{r_1}{2}\| \langle F, \Phi \rangle_{g_1} -  \langle G, \Psi \rangle_{w^2 g_2} \|_F^2
+ \frac{r_2}{2} \sum_{i=1}^N \int_{M_2}  \|\nabla_{\M_2} \psi_i \|^2 \mathrm{d}\M_2
+\frac{r_3}{2} \int_{\M_2} ||\nabla_{\M_2} w ||^2 \mathrm{d}\M_2, \\ 
&\hspace{3cm}\text{s.t.}\quad \int_{\M_2} \psi_i \psi_j \ w^2 \mathrm{d}\M_2  =\delta_{ij} \quad 
\text{and} \quad Area(\M_1)_{g_1} = Area(\M_2)_{w^2 g_2}  
\end{split}
\end{equation} 
where $Area(\M_1)_{\g_1} = \int_{\M_1} 1 \mathrm{d}\M_1$ and $Area(\M_2)_{\g_2} = \int_{\M_2} 1 w^2 \mathrm{d} \M_2$

\section{Discretization and Numerical Algorithms}
\label{sec:Algs}
In this section, we describe a discretization of the proposed variational model \eqref{eqn:finalmodel} using on triangular representation of surfaces. After that, we design a numerical algorithm to solve the proposed model based on proximal alternating minimization method. 

\subsection{Discretization of the Model}

The main method we use to discretize surfaces and differential operators is based on a finite element scheme similar to that developed in~\cite{reuter2006laplace,sun2009concise,dziuk2013finite}. Let $\{p_i\}_{i = 1}^{n}$ be a set of vertices sampled on the manifold $\M$. A surface can be discretized as a triple $\{P,E,T\}$ made of vertices ($P$), connected by edges ($E$) which form triangular faces $(T)$. We define the first ring of $p_i$, the set of all triangles which contain $p_i$ as $N(p_i)$. For each edge $E_{ij}$ connecting points $p_i$ and $p_j$, we define the angles opposite $E_{ij}$ as angles $\alpha_{ij}$ and $\beta_{ij}$.

We define a diagonal mass matrix, $\mathbf{M}$, a $n \times n$ positive definite matrix with entries given by:
$$\mathbf{M}_{ii} = \frac{1}{3} \sum_{\tau \in N(p_i)}  Area(\tau)$$
We use this simplified version, rather than the standard finite element discretization, for convince in order to avoid expensive factorizations later in our algorithm. We remark that the standard version can also be used in our algorithm at the cost of speed. The surface area can be approximated as $Area(\M) \approx \sum_{i=1}^n \textbf{M}_{ii}$. Similarly, a function $f:\M \rightarrow \RR$ with discretization $F:P \rightarrow \RR$, then we have the approximation $\int_{\M} f(x) \ d\M \approx 1^T\textbf{M}F = \sum_{i=1}^n f_i M_{ii}$. 
The stiffness matrix, $\textbf{S}$, is a $n \times n$ symmetric positive semidefinite matrix given by:
$$\mathbf{S}_{ij} = \sum_T \int_T \nabla_T e_i \cdot \nabla_T e_j = - \frac{1}{2}[\cot \alpha_{ij} (p_i) + \cot \beta_{ij} (p_i)] $$
where $e_i$ is a linear pyramid function which is $1$ at $p_i$ and zero elsewhere. The mass and stiffness matrices can be used to approximate the LB eigenvalue problem as: $\mathbf{S} f  =  \lambda \mathbf{M} f$ \cite{meyer2002discrete}. 

We remark that one can also work with point clouds representation instead of triangulated meshes. These definitions for the stiffness and mass matrices can be approximated by the point clouds method discussed in \cite{lai2013local}. The only change we would need to make is to use only the diagonal entries of the version of the mass matrix $\mathbf{M}$ proposed in their paper to populate the strictly diagonal version we employ here.

Suppose two surfaces $(\M_1,g_1),( \M_2,g_2) $ are represented by triangular meshes with the same number of points\footnote{In fact, we do not need to require that the surfaces have the same number of points, but doing so for now will allow for more convenient notation.}. We denote $\mathbf{M}^1, \mathbf{S}^1 \in \RR^{n\times n}$ as the mass and stiffness matrices of $\M_1$ and let $\Phi\in\RR^{n\times k}$ be the first $k$ LB eigenfunctions of $\M_1$, and $F \in\RR^{n\times \ell}$ be $\ell$ feature functions. 
Similarly, we write $\mathbf{M}^2, \mathbf{S}^2$ as the mass and stiffness matrices of  $\M_2$, $\Psi$ as the first $k$ LB eigenfunctions of $\M_2$ (under $g_2$), and $G$ as $\ell$ corresponding feature functions, ordered the same as in $F$. We also write $w^2$ as the discretized conformal factor on $\M_2$ and $\mathrm{diag}(w)$ as a diagonal matrix. 

Therefore, the discretized optimization model \eqref{eqn:map2} can be written as: 
\begin{equation}
\label{eqn:map2_discretization1}
\begin{split}
(w^*, \Psi^*) &= \arg\min_{w, \Psi} 
\frac{r_1}{2} \| F^T \mathbf{M}^1 \Phi - G^T  \mathrm{diag}(w) \mathbf{M}^2 \mathrm{diag}(w) \Psi \|_F^2 
 + \frac{r_2}{2} \tr (\Psi^T \mathbf{S}^2 \Psi) 
+ \frac{r_3}{2}w^T \textbf{S}_2 w, \quad 
\\ 
&\hspace{3cm}\text{s.t.} \quad \Psi^T   \mathrm{diag}(w) \mathbf{M}^2 \mathrm{diag}(w) \Psi = \textbf{I}_k, 
\quad \text{and} \quad w^T \textbf{M}_2 w = A 
\end{split}
\end{equation}
Here $\textbf{I}_k$ is the $k \times k$ identity matrix and $A = \sum_{i=1}^n \textbf{M}_{1}(i,i) $. Since $\textbf{M}_2$ is symmetric positive definite and diagonal, we can easily calculate the matrix decomposition $\textbf{M}_2 = \textbf{L}^T \textbf{L}$. If we also substitute $\bar\Psi = \textbf{L}~ \mathrm{diag}(w) \Psi $, then \eqref{eqn:map2_discretization1} can be written as:
\begin{equation}
\label{FinalDisc}
\begin{split}
(w^*, \bar \Psi^*)& = \arg\min_{w, \bar \Psi}  \E(w, \bar \Psi ) = 
\frac{r_1}{2} \| F^T \textbf{M}_2 \Phi - G^T  \mathrm{diag}(w) \textbf{L}^T \bar\Psi \|_F^2 
 + \frac{r_2}{2} \tr (\bar\Psi^T \bar{\textbf{S}}^2(w) \bar\Psi)
 + \frac{r_3}{2} w^T \textbf{S}_2 w, \quad
\\
&\hspace{4cm}  \text{s.t.} \quad \bar \Psi^T \bar \Psi = \textbf{I}_k 
\quad   \text{and} \quad
 w^T \textbf{M}_2 w = A 
\end{split} 
\end{equation}
where $ \bar{\textbf{S}}^2(w) = (\textbf{L}^T)^{-1}\mathrm{diag}(w)^{-1} \textbf{S}_2 \mathrm{diag}(w)^{-1} \textbf{L}^{-1}.$ Note that this parameterization of the problem moves the conformal factor $w$ out of the orthogonality constraint (and into $\bar{\textbf{S}}$). We will soon see that, for any fixed $\bar \Psi$, this will make the problem for $w$ easier to solve.

\subsection{Numerical Optimization}
\label{subsec:optimization}
The two variables $w$ and $\bar\Psi$ in \eqref{FinalDisc} make the optimization problem different from orthogonality constrained problems solved by nonconvex alternating direction method of multipliers (ADMM) methods considered in \cite{lai2014splitting,chenAugmented,wang2015global,kovnatsky2016madmm}. Rather than solve this problem directly for $\bar \Psi$ and $w$ simultaneously by directly minimizing \eqref{FinalDisc}, we employ a method based on the framework of proximal alternating minimization (PAM) method~\cite{attouch2010proximal}. The idea of PAM is to update variables using a Gauss-Seidel method with proximal regularizations in the following way.

Let $\S = \{\bar \Psi \in\RR^{n\times k} ~|~ \bar \Psi^T \bar \Psi = \textbf{I}_k \}$ and $\mathcal{W} = \{w\in\RR^n~|~w^T \textbf{M}_2 w = A \}$. We also define indicator functions
$$\delta_\S(x) = \left\{\begin{array}{cc}0,& \text{if } x\in\S\\ +\infty, &\text{otherwise}\end{array}\right.,\qquad  \delta_\W(x) = \left\{\begin{array}{cc}0,& \text{if } x\in\W\\ +\infty, & \text{otherwise}\end{array}\right. .$$Then it is clear that $\delta_\S$ and $\delta_\W$ are semi-algebraic functions as $\S$ and $\W$ are zero sets of polynomial functions~\cite{attouch2013convergence}. Therefore, we write an equivalent form of \eqref{FinalDisc} as
\begin{equation}
\label{eqn:FinalDisc_PAM}
(w^*, \bar \Psi^*) = \arg\min_{w, \bar \Psi} \E(w, \bar \Psi ) + \delta_\S(\bar\Psi) + \delta_\W(w).
\end{equation}
Using the PAM method, we have the following iterative scheme
\begin{eqnarray} 
\label{optstep}
\left\{
\begin{aligned}
\bar \Psi^{j+1} &= \arg\min_{\bar \Psi}\E(w^j, \bar \Psi )  + \frac{1}{2\eta} ||\bar \Psi - \bar \Psi^j ||^2,  \quad \text{s.t.} \quad \bar \Psi^T \bar \Psi = \textbf{I}_k \\ 
w^{j+1} &= \arg\min_{w} \E(w, \bar \Psi^{j+1} )  + \frac{1}{2\eta} ||w - w^j ||^2,  \quad \text{s.t.} \quad w^T \textbf{M}_2 w = A 
\end{aligned}\right.
\end{eqnarray}
Here $\eta$ is a step size parameter. 
Essentially, these proximal terms penalizes large step sizes in and prevents the algorithm from ``jumping" between multiple local minimums. The addition of these proximity terms allows us to analyze our proposed method in the framework of the PAM algorithm~\cite{attouch2010proximal}. It has been shown in \cite{attouch2010proximal, attouch2013convergence,bolte2014proximal} that such proximal terms can guarantee the solutions generated at each step converge to a critical point o          bbbbbbbbf the objective function. Formally, we have the following convergence theorem in accordace with Theorem 9 in~\cite{attouch2010proximal}.
\begin{theorem}\label{thm:convergence}
Let $\{w^j,\bar\Psi^j\}$ be the sequence  produced by \eqref{optstep}, then the following statement hold:
\begin{enumerate}
\item $\displaystyle \E(w^{j+1},\bar\Psi^{j+1}) + \frac{1}{2\eta} ||\bar \Psi^{j+1} - \bar \Psi^{j} ||^2 +  \frac{1}{2\eta} ||w^{j+1} - w^{j} ||^2 \leq \E(w^{j},\bar\Psi^{j}), ~\forall j\geq 0$.
\item $\displaystyle \sum_{j=1}^{\infty} (\|w^{j} - w^{j-1}\|^2 + \|\bar\Psi^{j} - \bar\Psi^{j-1}\|^2) < \infty$.
\item $\{w^j,\bar\Psi^j\}$ converges to a critical point of $\E(w,\bar\Psi)$.
\end{enumerate}
\end{theorem}

\begin{proof}
To prove this, we show that our model obeys the conditions required for local convergence of PAM in \cite{attouch2010proximal}. To do so, we need: 

(1) Each of term in the objective which contains only one primal variable is bounded below and lower semicontinous.

(2) Each term in the objective which contains both variables is $C^1$ and has a locally Lipschitz continuous gradient.

(3) The objective satisfies the Kurdyka-Lojasiewicz (KL) property. \\
It is immediately clear that that the first two properties are satisfied by our objective. Furthermore, it is known that all semi-algebraic functions have KL property \cite{attouch2010proximal,attouch2013convergence, chenAugmented}. Our objective is semi-algebraic so we can guarantee local convergence of the proposed optimization method. 
\end{proof}

We use the augmented Lagrangian method to solve the constrained sub-optimization problem for $w$ in \eqref{optstep}. For convenience, let's write
\begin{equation}\label{Lap}
\mathcal{L}(\bar \Psi,w;b) = 
\E(w, \bar \Psi )
+ \frac{r_4}{2} \Big(  w^T \textbf{M}_2 w - A   + b\Big )^2
\end{equation}
Overall, we solve \eqref{FinalDisc} in the following way by hybridizing PAM with the augmented Lagrangian method.
\begin{eqnarray} \label{optsteps}
\left\{\begin{split}
 \bar \Psi^{j+1} &= \arg\min_{\bar \Psi} \mathcal{E} (w^j, \bar \Psi) + \frac{1}{2\eta} ||\bar \Psi - \bar \Psi^j ||^2  \quad \text{s.t.} \quad \bar \Psi^T \bar \Psi = \textbf{I}_k\\ 
w^{j+1} & \leftarrow  \begin{cases}
 \displaystyle w^{j+1,s+1} = \arg \min_{w}  \mathcal{L} (w, \bar \Psi^{j+1};b^{j+1,s}) + \frac{1}{2\eta} ||w - w^{j} ||^2  \\ 
 \displaystyle  b^{j+1,s+1} =  b^{j+1,s} +  (\w^{j+1,s+1})^T \textbf{M}_2 w^{j+1,s+1} -  A  .\\
\end{cases}
\end{split}\right.
\end{eqnarray}

The subproblem for minimizing $\bar\Psi$ requires a different approach to be solved. The main challenge in minimizing the first sub-optimization problem is the nonconvex orthogonality constraints. Recently, several approaches have been developed to solve orthogonally constrained problems in feasible or infeasible ways~\cite{wen2013feasible,lai2014splitting,chenAugmented,wang2015global,kovnatsky2016madmm}. For our implementation, we have chosen the feasible approach developed in \cite{wen2013feasible} which uses a curvilinear method based on the Cayley transform together with Barzilai-Bowein step size line search. This method updates variables along a geodesic curve on the Stiefel manifold, a geometric description of the orthogonality. It preserves the orthogonality constraints and guarantees convergence to critical points in our scenario. More precisely, given a feasible starting point $\bar \Psi^s$  and the coordinate gradient $Y^s$ at this point, the update scheme is as follows:
\begin{equation}
\label{eqn:curvilnear}
\left\{
\begin{split}
D^s &= Y^s (\bar \Psi^s)^T - \bar\Psi^s(Y^s)^T \\ 
Q^s &= (I+\frac{dt}{2}D^s)^{-1}(I-\frac{dt}{2} D^s) \\ 
\bar \Psi^{s+1} &= Q^s \bar \Psi^s 
\end{split}\right.
\end{equation}
Here $dt$ is a step size parameter chosen by the Barzilai-Bowein criteria developed in \cite{barzilai1988two}.  Although convergence to a global minimum is not guaranteed, this method has proven effective for our purposes and only requires the computation of the objective function and its coordinate gradient $Y^s$ with respect to $\bar\Psi$ at each step provided by:

\begin{eqnarray} \label{grad Epsi}
\begin{aligned}
\nabla_{\bar\Psi} \left(\E(w,\bar\Psi) +\frac{1}{2\eta} ||\bar \Psi - \bar \Psi^j ||^2 \right) = & - r_1 G^T \mathrm{diag}(w) \textbf{L}^T \Big(F^T \textbf{M}_1\Phi - G^T \mathrm{diag}(w) \textbf{L}^T\bar \Psi\Big)\\
& + r_2 \bar{\textbf{S}}^2 \bar \Psi + \frac{1}{\eta}(\Psi - \bar\Psi^j)
\end{aligned}
\end{eqnarray}

The subproblem for $w$ (as written in \eqref{optsteps}), on the other hand is smooth and unconstrained.
For our implementation, we use the well known quasi-Newton BFGS algorithm \cite{bazaraa2013nonlinear}. The gradient of objective function with respect to $w$ can be written as:
\begin{eqnarray} \label{grad Ew}
\begin{aligned}
\nabla_{w}\left(\mathcal{L}(w,\bar\Psi;b) + \frac{1}{2\eta} ||w - w^{j} ||^2 \right) &= r_1 \ diag \Big(G^T(F^T\textbf{M}_1 \Phi-G w \textbf{L}^T \bar\Psi)) \bar\Psi^T \textbf{L} \Big) \\
& + r_2 \ diag \Big(\Psi \Psi^T \textbf{S} w^{-1} \Big)\odot w^{-2} \\
 & + r_3 \textbf{S}_2 w 
 + r_4 \Big( w^T \textbf{M}_2 w -  A + b \Big)  \textbf{M}_2 w + \frac{1}{\eta}(w - w^{j})
\end{aligned}
\end{eqnarray}
where $diag \big(\cdot \big)$ denotes the diagonal of the matrix, $\odot$ signifies element-wise Hadamard product and $w^{-2}$ is the inverse of diagonal matrix $w$ multiplied with itself. 

\subsection{Computation of Point-to-Point Map}

One naive way to compute a point-to-point map is to find the functional map by using the final deformed manifold and its LB eigensystem with respect to the deforamtion. However, this may not work well because of the ambiguity of LB eigensystem. Additional effort is needed to handle possible ambiguity of LB eigensystem such as the method discussed in~\cite{lai2017multi}. As an advantage of the proposed method, the resulting basis generated by the proposed algorithm (recovered as $\Psi^* = A^{-1} w \bar \Psi$) to will naturally correct ambiguities of LB eigensystem. This is similar to the method discussed in~\cite{kovnatsky2013coupled}.  Thus, we can compute the functional map as $F_T(h) = \sum_{i=1}^k  (\int_{\M_1} h \phi_i \ \mathrm{d}\M_1)  \psi^T_i = \Psi \Phi^T \textbf{M}_1h$. However, this method is still quite inefficient and may be sensitive to small errors in the resulting basis. 

Instead, after we recover the final basis from our method, we can compute the point-to-point map between the two surfaces by comparing the values of each of the basis functions. This is essentially the same scheme presented in \cite{ovsjanikov2012functional}, but applied to our new basis. We use a KNN search (with $K = 1$) to match rows of $\Phi$ and $\Psi$. This requires a search of $n$ points in $k$ dimension, but is much more efficient and accurate than using the delta function approach described in the previous paragraph. Other methods used to refine functional maps such as \cite{rodola2015point} can be applied in this setting without changes.

We summarize our numerical method for nonisometric surface registration as Algorithm \ref{alg:LBBasisPursuit1}.

\begin{algorithm2e}[h]
\caption{LB Basis Pursuit Algorithm}
\label{alg:LBBasisPursuit1}
\SetKwInOut{Input}{input}\SetKwInOut{Output}{output}
\SetKwComment{Comment}{}{}

\KwIn{Triangulated surfaces $\M_1$ and $\M_2$ and list of known corresponding functions $F$ and $G$.}

\KwOut{$\Psi^*$, $w$, point-to-point map}

Compute stiffness and mass matrices for each surface: $\textbf{M}_1, \textbf{M}_2, \textbf{S}_1, \textbf{S}_2$\;
Use stiffness and mass to calculate LBO eigensystems: $\textbf{M}_1 \Phi = \lambda \textbf{S}_1 \Phi$\;
Initialize: Let $\Psi^0$ be the LB eigenfunctions of target surface: $\textbf{M}_2 \Psi = \lambda \textbf{S}_2 \Psi$\;
Compute $\bar\Psi^0 = \textbf{L} w \Psi$\;
\While{not converged}{
Update $\displaystyle \bar\Psi^{j+1}= \arg\min_{\bar\Psi}\mathcal{E} (w^j, \bar \Psi) + \frac{1}{2\eta} ||\bar \Psi - \bar \Psi^j ||^2$ using the curvilinear search algorithm \eqref{eqn:curvilnear}\;
\While{ $s \leq \ell$ }{
Update $\displaystyle w^{j+1,s} = \arg \min_{w} \mathcal{L} (w, \bar \Psi^{j+1};b^{j+1,s}) + \frac{1}{2\eta} ||w - w^{j} ||^2$ using BFGS\;
$\displaystyle b^{j+1,s+1} =  b^{j+1,s} + (w^{j+1})^T \textbf{M}_2  w^{j+1} - A$\;
}
$w^{j+1} = w^{j+1,s}$\;
 }
Recover $\Psi^* = w \textbf{L}^{-1} \bar \Psi $\;
Compute correspondence map with KNN-search of coefficient space 
\end{algorithm2e}

\section{Discussion}
\label{sec:Discussion}
In this section, we discus our choice of feature functions, as well as ways to overcome problems which may arise from the non-convexity of the proposed optimization problem. In addition, we present a novel way to jointly measure the quality of the correspondence and alignment of the bases without any prior knowledge about the ground truth of the point-to-point map.   

\subsection{Choice of Feature Functions}
The simplest, and in many applications, most natural features to choose for $F$ and $G$ are indicator functions. Let  $\{\chi^{1}_i\}_{i=1}^k$ be a set of points on $\M_1$ and $\{\chi^{2}_i \}_{i=1}^k$  be a corresponding set on  $\M_2$. We can view each $f_i$ and $g_i$ as a $\delta$-function on $\M_1$ and $\M_2$ respectively to indicate these landmarks.

Another option is to use heat diffusion functions. Given a corresponding pair of points we can use delta functions to define an initial condition and solve the heat diffusion problem $\frac{\partial u}{\partial t}(\x) = \Delta u(\x,t)$ using the Crank-Nicholson scheme 
$ \Big( \textbf{M}+\frac{dt}{2} \textbf{S} \Big)u^{i+1} = \Big( \textbf{M} -\frac{dt}{2}\textbf{S} \Big) u^{i}$
where $dt$ is a step size parameter.
By taking ``snap shots'' (solutions of the equation for various $t$ values) of $u$ at different time values we can generate multiple functions from a single corresponding pair. This choice allows for a multi-scale selection of features. 

The wave kernel signature (WKS) has also been used for characterizing points on non-rigid three dimensional shapes \cite{aubry2011wave}. These functions are defined as the solutions to the Schrodinger equation: $\frac{\partial u}{\partial t}(\x) = i \Delta u(\x,t)$  at different points on the surface. Given two corresponding points we can solve the equation at each point and use these as our corresponding functions. However, the solutions to these equations are highly dependent on both local and global geometries of the manifold. Because of this, they are only suitable for shape correspondence when the shapes are very similar and, in general, do not work well for non-nearly-isometric problems. The same problem exists for heat diffusion features, however, in general heat diffusion tends to be much more stable with respect to local deformations.%

\subsection{Reinitialization Schemes}
Although we have shown that the proposed PAM based optimization algorithm converges to a critical point of the objective function, it is still challenge to achieve global optimizers as the problem is non-convex. In practice, we have found that the numerical results can often be improved in terms of both accuracy and speed of computation by adding a simple reinitialization scheme to our algorithm. The motivation for the scheme comes from an observation that if we know the exact conformal deformation $w^2$ and the source surface has a simple eigensystem (no repeated eigenvalues), then the LB eigensystem of $(\M_1, g_1)$ is the same as the LB eigensystem of $(\M_2,w^2 g_2)$ up to a change in sign.  
With this in mind, we propose to reinitialize the $\Psi$ problem by resetting $\Psi$ to be the solution to weighed eigenproblem $\textbf{S}_2 \Psi = \Lambda  \mathrm{diag}(w^2) \textbf{M}_2 \Psi$ . Computationally, to avoid introducing ambiguities of LB eigensystem by calling standard eigen-solvers, we solve a discrete counterpart to \eqref{defeig} as
$\displaystyle \min_{\Psi} ~ \tr (\bar\Psi^T \bar{\textbf{S}}^2(w) \bar\Psi) , 
~ \text{s.t.} ~ \bar \Psi^T  \bar \Psi = \textbf{I}$
based on the curvilinear search method discussed in Section \ref{subsec:optimization} with the current eigensystem, $\bar \Psi^{j+1}$, as an initial guess for this problem. By using $\bar \Psi^{j+1}$ as warm start for the eigenproblem we can avoid re-introducing sign or multiplicity ambiguities into the problem which our algorithm has already resolved. 

If using heat diffusion, wavelet kernel signatures, or any other functions which are defined based on local geometry as the  input feature functions, then we also need to recalculate these functions with respect to the conformally deformed metric. For example, if we are using heat diffusions, we can recompute the heat diffusion functions on the deformed manifold $(\M_2,w^2 g_{2})$ by multiplying the mass matrix by $w^2$ in the Crank-Nicholson scheme:
$ \Big(\textbf{M}_2 \mathrm{diag}(w^2)+\frac{dt}{2}\textbf{S}_2 \Big)u^{i+1}
 = \Big(\textbf{M}_2 \mathrm{diag}(w^2)-\frac{dt}{2}\textbf{S}_2\Big)u^{i}$,  
A similar re-computation technique can be applied if using WKS features.

With this re-initialization procedure, we propose a modified version of our numerical solver as Algorithm~\ref{alg:LBBasisPursuit2}. 

\begin{algorithm2e}
\caption{Basis Pursuit Algorithm with reinitialization}
\label{alg:LBBasisPursuit2}
\KwIn{Set of vertices and faces of source ($\M_1$) and target ($\M_2$) manifolds and list of known corresponding functions $F$ and $G$}
\KwOut{$\Psi^*$, $w^*$, point-to-point correspondence map}
Compute stiffness and mass matrices for each surface: $\textbf{M}_1, \textbf{M}_2, \textbf{S}_1, \textbf{S}_2$\;
Use stiffness and mass to calculate LBO eigensystems: $\textbf{M}_1 \Phi = \lambda \textbf{S}_1 \Phi$\;
Compute corresponding feature functions $F$ and $G$ on $\M^1$ and $\M^2$ respectively\;
Initialize: Let $\Psi^0$ be the LBO eigenfunctions of target surface: $\textbf{M}_2 \Psi = \lambda \textbf{S}_2\Psi$\;
\While{number of re-initialization steps complete $<$ max number of re-initializations}{
Update $\displaystyle \bar\Psi^{j+1}= \arg\min_{\bar\Psi}\mathcal{E} (w^j, \bar \Psi) + \frac{1}{2\eta} ||\bar \Psi - \bar \Psi^j ||^2$ using the curvilinear search algorithm \eqref{eqn:curvilnear}\;
\While{ $s \leq \ell$ }{
Update $\displaystyle w^{j+1,s} = \arg \min_{w} \mathcal{L} (w, \bar \Psi^{j+1};b^{j+1,s}) + \frac{1}{2\eta} ||w - w^{j} ||^2$ using BFGS\;
$\displaystyle b^{j+1,s+1} =  b^{j+1,s} +  (w^{j+1})^T\textbf{M}_2 w^{j+1}  - A$\;
}
$w^{j+1} = w^{j+1,l}$\;

	\If {update $<$ tolerance}{
		Re-Initialize $\bar \Psi$ as $\displaystyle \argmin_{\bar\Psi} ~ \tr (\bar\Psi^T \bar{\textbf{S}}^2(w^{j+1}) \bar\Psi) , ~ \text{s.t.} ~ \bar \Psi^T \textbf{M}_2 \bar \Psi = \textbf{I}$\;
		\If {Using feature functions which depend on local geometry}
		{Re-Compute features using $\textbf{M}_2 \mathrm{diag}(w^2)$ as Mass matrix}
	}
}
Compute correspondence map with KNN-search of coefficient space;

\end{algorithm2e}

\section{Numerical Experiments}
\label{sec:Results}
In this section, we apply our algorithm to several problems. We begin by working on a typical non-isomorphic matching problem for a pair of shapes with a large deformation: a horse and an elephant. We preform tests showing the effectiveness of our approach given different amounts of landmark points, and demonstrate robustness with respect to noise both on the manifold and in the initial correspondences. Finally, we conduct experiments on the Faust benchmark data set \cite{bogo2014faust}. All numerical experiments are implemented in MATLAB on a PC with a 32GB RAM and two 2.6GHz CPUs.

In all of our experiments, we use randomly chosen correspondence points to create indicator functions as the input features.  The first 100 non-trivial LB eigenfunctions are chosen to calculate the coefficient matching term, as well as for computing the final correspondence. We set $r_1=10,\ r_2=10,\ r_3=1,\ r_4=.01, \ell = 1$ for all experiments, even though the data sets and experimental conditions are very different. This choice of $r_1$ and $r_2$ allows the coefficient matching terms and eigenfunction term to balance each other out, with the choice of $r_3$ still being large enough to preserve the area constraint. $r_4$ is chosen to be small so that the harmonic energy, which tends to be quite large, does not dominate the others. We observe that our algorithm is robust to different choices of the parameters. 

\subsection{A large deformation pair: Horse to Elephant}
\begin{figure}[htp]
\begin{minipage}{0.55\linewidth}
\begin{center}
\includegraphics[width=.9\linewidth]{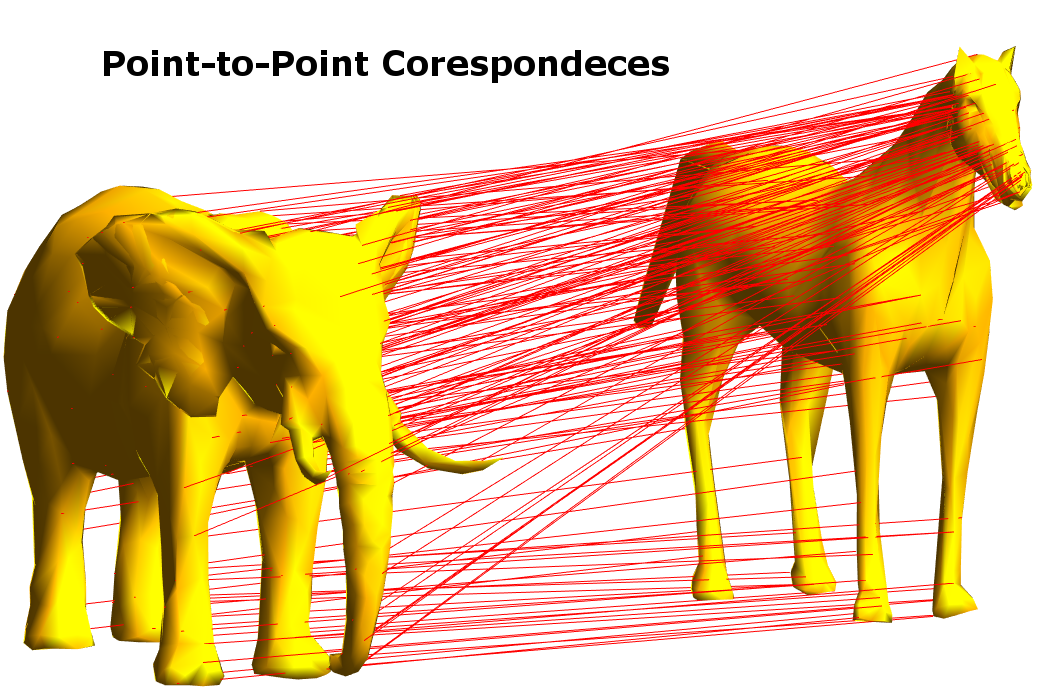}
\end{center}
\end{minipage}
\begin{minipage}{0.49\linewidth}
\begin{center}
\includegraphics[width=.95\linewidth]{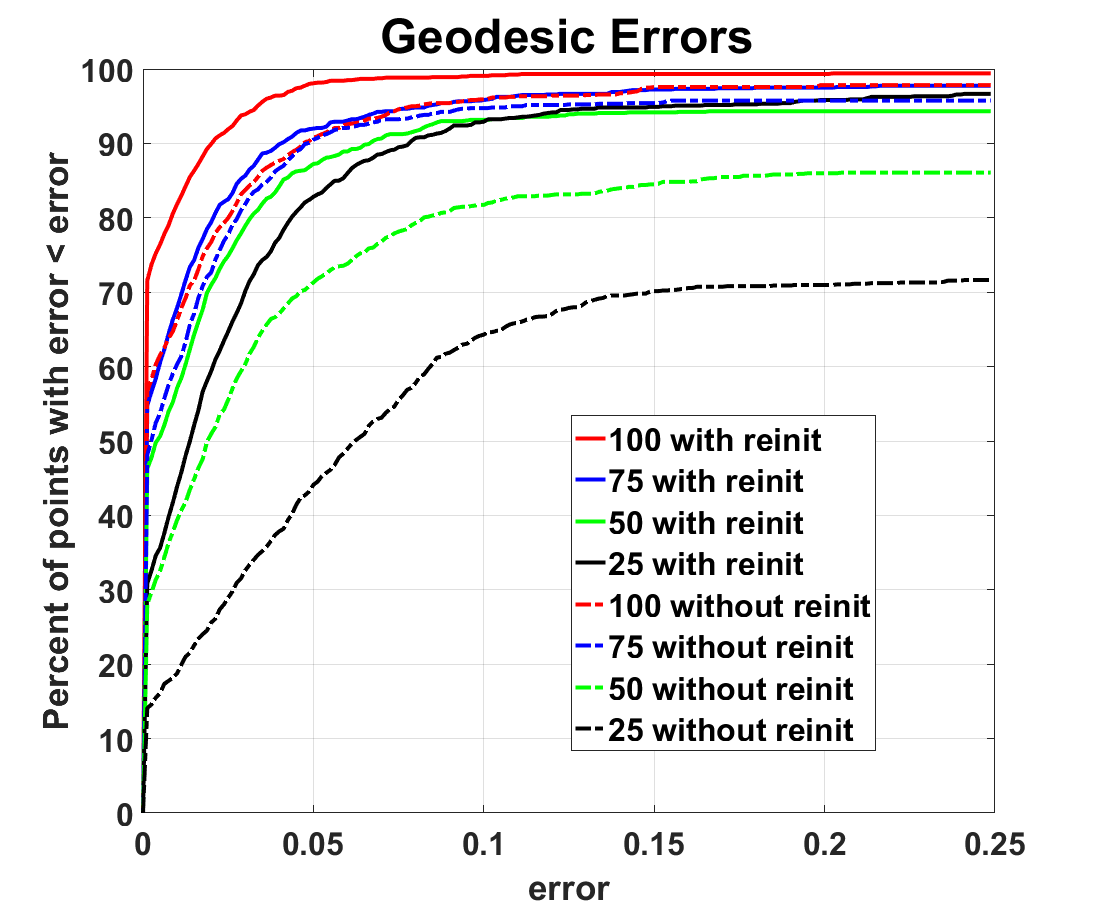}
\end{center}
\end{minipage}
\caption{Left: visualization of point-to-point map. Right: normalized geodesic errors for various numbers of landmarks with and without reinization.}
\label{P2Ppic}
\end{figure}
The first experiment is designed to test the effectiveness of the proposed method on a pair of shapes with large deformation. Each surface, a horse and an elephant, is represented by a mesh with 1200 points. One of the challenges in this pair is the large deformations in the sharp corner and elongated regions including ears, teeth, noses and tails on the horse and elephant surfaces. Those regions make the registration problem very difficult.To demonstrate the efficacy of our approach, we perform this experiment under several different conditions. Our algorithm produces excellent results given a sufficient number of landmarks, and it still finds reliable correspondences given limited landmarks. We also show that using our reinitialization scheme (Algorithm \ref{alg:LBBasisPursuit2})  produces a more accurate map than without this extra step (Algorithm \ref{alg:LBBasisPursuit1}). 

\begin{figure}[t]
\begin{minipage}{.46\linewidth}
\centering
\includegraphics[width=.9\linewidth]{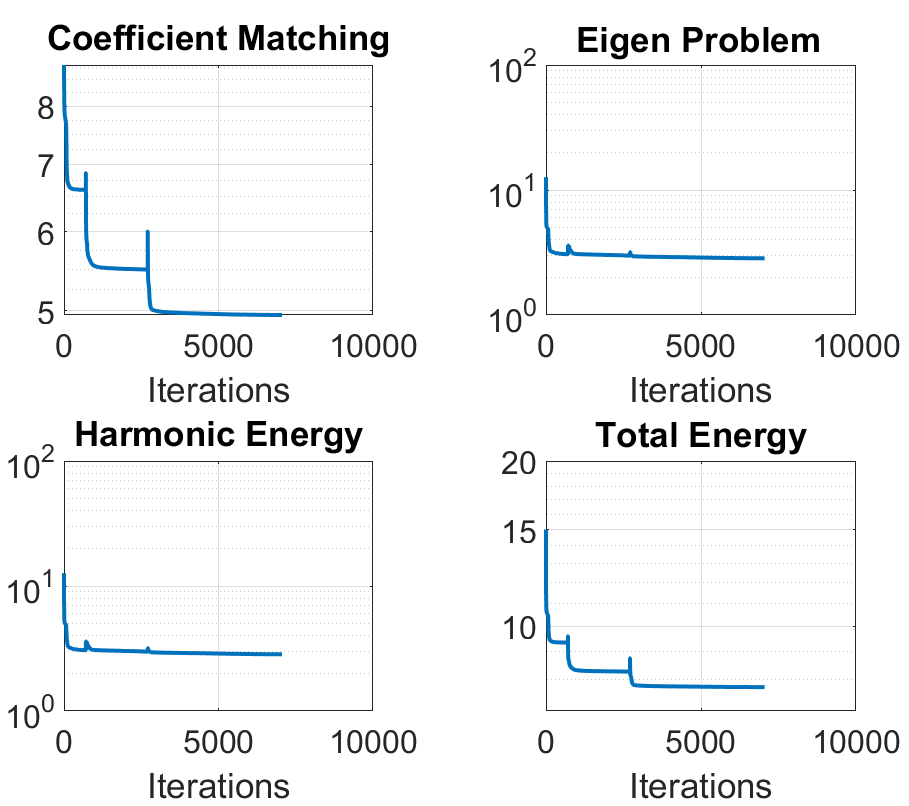}
\end{minipage}\hfill
\begin{minipage}{.46\linewidth}
\centering
\includegraphics[width=.9\linewidth]{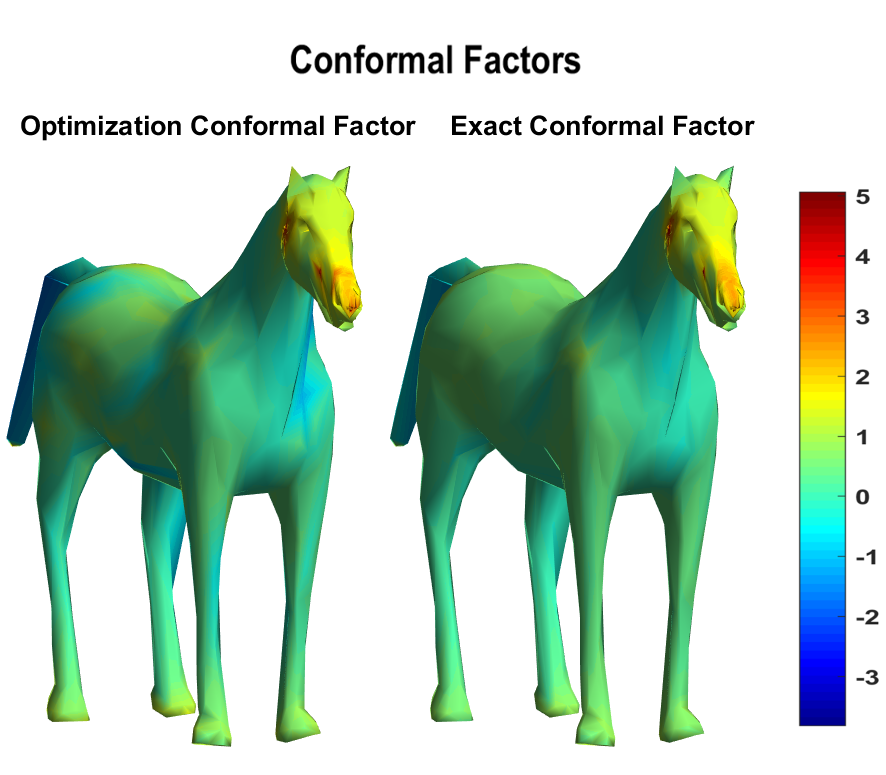}
\end{minipage}
\caption{Left: convergence curves of our method. The coefficient matching term measures: $\| F^T \textbf{M}_1 \Phi - G^T  \textrm{diag}(w) \textbf{L}^T \bar\Psi \|_F$. The eigen problem is: $(\Psi^T \textbf{S}_2 \Psi)$ and the harmonic energy measures: $w^T \textbf{S}_2 w$ and the total energy is the entire model derived in \eqref{FinalDisc}.  Right: resulting and exact conformal factors.}
\label{convergence}
\end{figure}
Figure \ref{P2Ppic} shows the results of using 100, 75, 50 and 25 known landmark points with and without our reinitialization scheme.  To qualitatively measure the mapping quality, we calculate the normalized geodesic distance from the point on the target surface produced by the map to ground truth following the Princeton Benchmark method \cite{kim2011blended}. These distances are collected into a cumulative error on the right of Figure \ref{P2Ppic} where the $y$-axis measures the percent of points whose distances are less than or equal to the $x$-axis value. For example, in the case of $100$ known landmarks, our algorithm matches over $70\%$ of the points to exact correct point and more than $95\%$ within a $5\%$ error margin. 
\begin{figure}[H]
\begin{center}
\includegraphics[width=.9\linewidth]{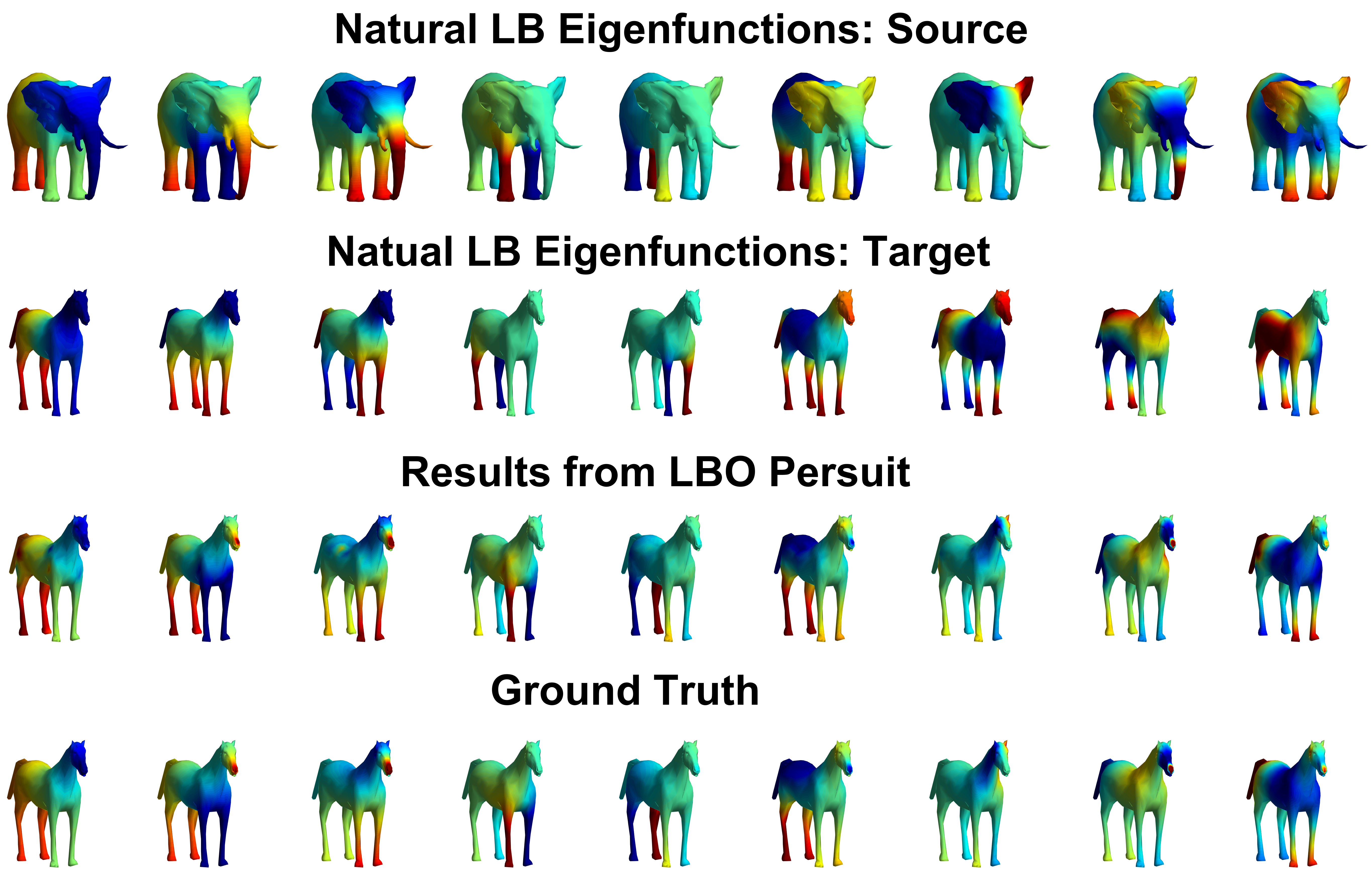}
\end{center}
\vspace{-0.5cm}
\caption{First two rows: The first 9 non-trivial natural LB eigenfunctions of manifolds.
Third row:  results from the proposed basis pursuit algorithm. 
Fourth row: ground truth.}
\label{H2Eeigs}
\end{figure}

The left of Figure \ref{convergence} shows the convergence of the objective function and illustrates the effectiveness of the reinitialization step. We plot the three terms in the objective function separately as well as the overall objective. We typically observe that the convergence curves in the coefficient matching and total energy flatten quickly as the algorithm tends to a local minimizer. However, each reinitialization significantly reduces the objective function, thus effectively overcoming the non-convexity of the model. We further demonstrate the validity of our algorithm by examining the resulting conformal factor. In the right image of Figure \ref{convergence}, we show the conformal factor calculated by our algorithm as well as the ground truth. The ground truth conformal factor is calculated by using the ground truth point-to-point map to compare the area of the first ring structure around each point on the source and target surface.  Here we plot $u$ where $w^2=e^{2u}$ for better visualization. From this figure we can confirm that the conformal mapping our algorithm produces is very close to the true factor. 

\begin{wrapfigure}{r}{0.6\textwidth} 
\centering
\includegraphics[width=1\linewidth]{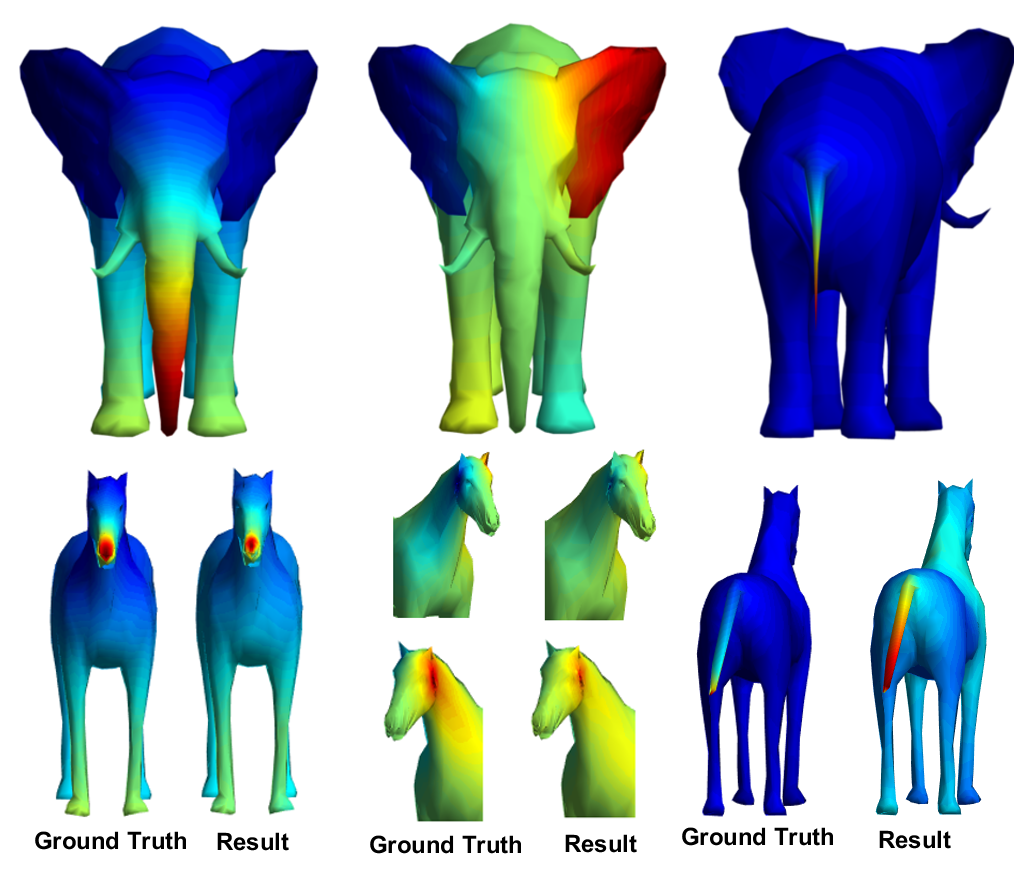}
\caption{Top row: 9th, 11th and 44th natural LB eigenfunctions on source. Bottom row: results and ground truth.}
\label{closeups}
\end{wrapfigure}

Since the elephant and horse are dramatically different shapes, the large dissimilarity of their natural LB eigenfunctions (first two rows of Figure \ref{H2Eeigs}) cannot be expected to produce meaningful correspondence. However, our model overcomes this by capturing the conformal deformation between the surfaces. As a result, the basis computed for the horse (target surface) by our model is consistent with the LB eigenfunctions of the elephant (source surface). We further compare these results with the ground truth which is calculated through the push forward of the LB eigenfunction of the source to the target surface using the \textit{a priori} map.  Figure \ref{closeups}  highlights the consistency of the produced bases on several highly distorted regions. Specifically, we focus on each of the ears, the nose/trunk and the tails. From Figures \ref{H2Eeigs} and \ref{closeups}, we can see that our approach produces a new basis on the target that aligns very closely to the natural LB basis on the sources manifold. This is the reason that accurate registration results can be obtained using the new basis.

\subsection{Noisy Data}
\begin{figure}[h]
\begin{minipage}{.49\linewidth}
\includegraphics[width=1\linewidth]{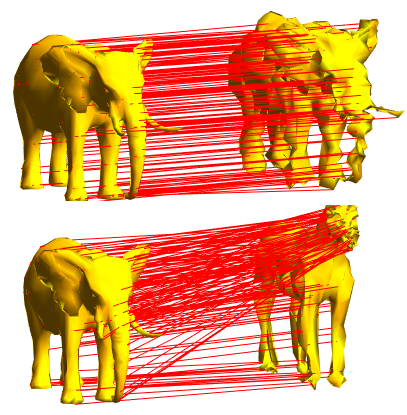}
\end{minipage}
\begin{minipage}{.49\linewidth}
\includegraphics[width=1\linewidth]{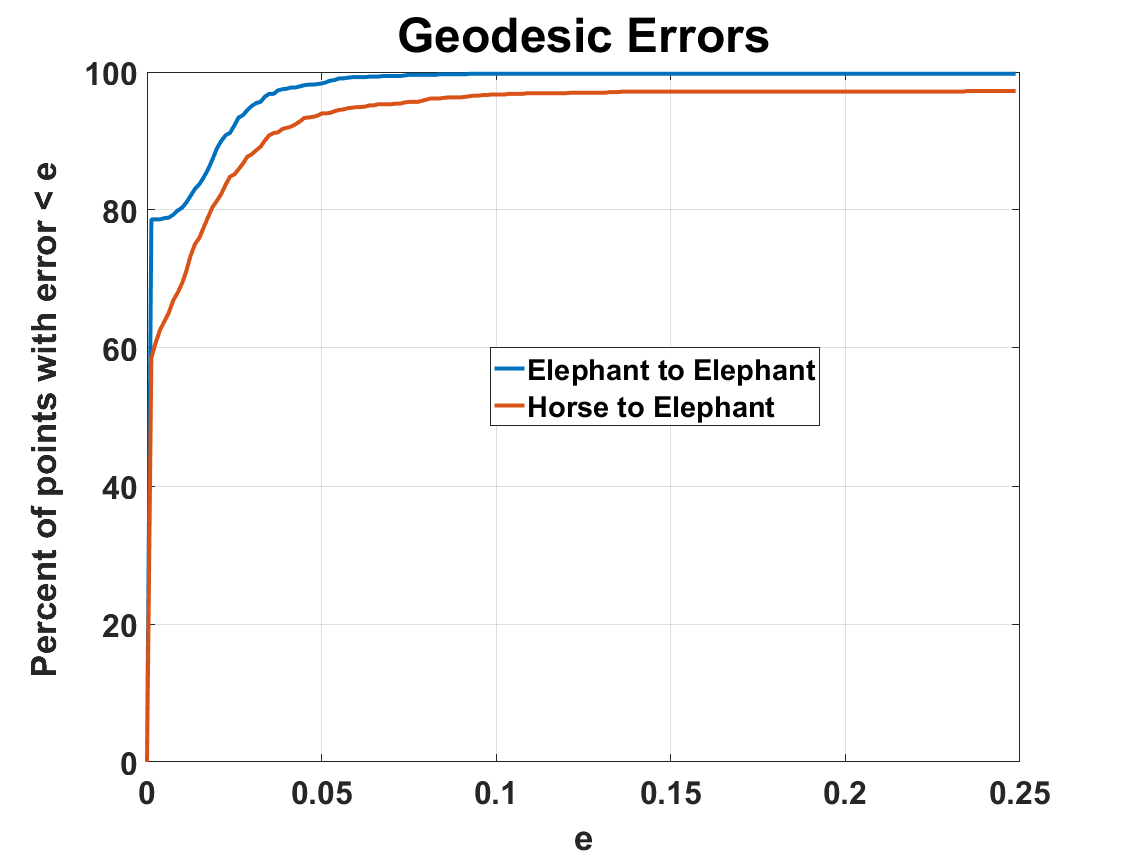}
\end{minipage}
\caption{Left: point-to-point maps for noisy data. Right: normalized geodesic errors for noisy data. }
\label{H2Enoise}
\end{figure}

In this experiment, we demonstrate that our algorithm can handle noisy data. Since noise on the surfaces can be viewed as local deformations, our algorithm is automatically robust to geometric noise. Medical scans often have noise resulting from the imaging instruments and segmentations. Our model can solve registration problems for this type of data. To demonstrate this, we generate noisy data by adding noise along the normal of each point. Figure \ref{H2Enoise} shows the results of two experiments: a noisy elephant to an elephant and a noisy horse to an elephant. We observe that our algorithm still produces very accurate results despite this noise.

\subsection{Perturbed Landmarks}
\begin{figure}[h]
\begin{minipage}{.49\linewidth}
\includegraphics[width=.99\linewidth]{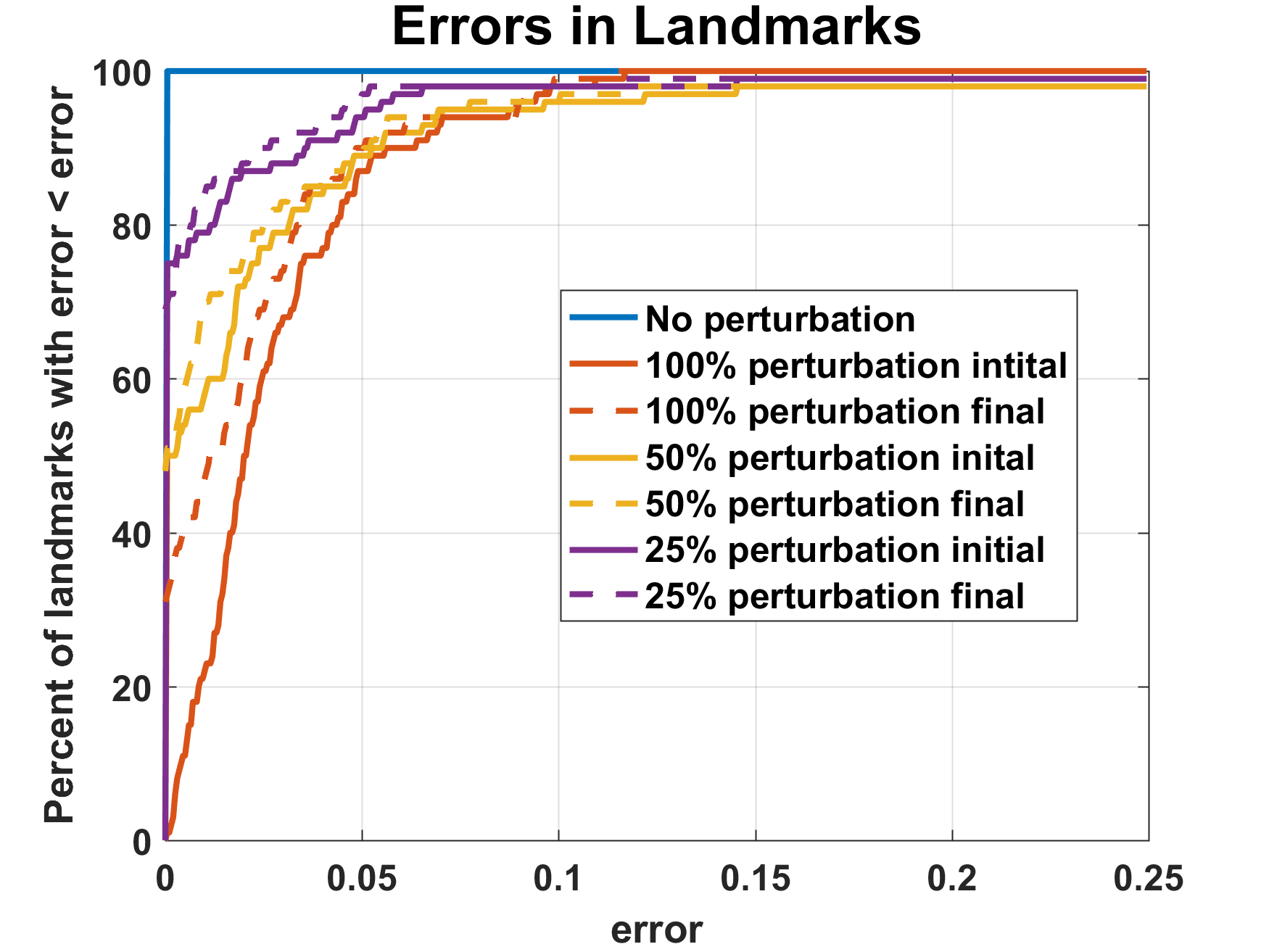}\\
\end{minipage}
\begin{minipage}{.49\linewidth}
\includegraphics[width=.93\linewidth]{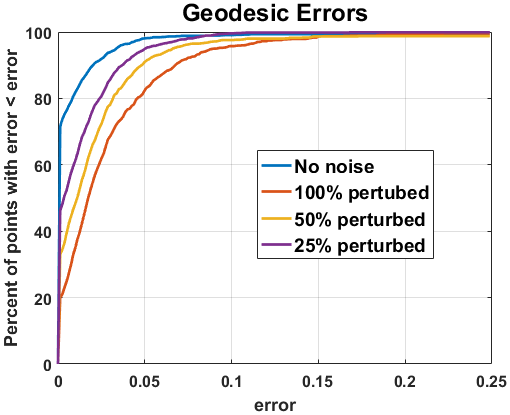}\\
\end{minipage}
\caption{Left: Initial perturbations to landmarks and final error of landmarks.  Right: Final registration geodesic errors for all points using perturbed landmarks.}
\label{noisycompgeo}
\end{figure}

We also demonstrate the robustness of our algorithm to landmark perturbations. Working again on the horse and elephant, we test cases where the landmarks are perturbed to another vertex within the first ring. The magnitude of these perturbations depends on the uniformality and meshing of the surface. The left graph in Figure \ref{noisycompgeo} shows the size of the perturbations of the landmarks points as well as the error in their final mapping. The right graph in Figure \ref{noisycompgeo} compares the geodesic error of the for all points when 25\%, 50\% and 100\% of the landmarks points are perturbed. From these tests we conclude that our method can successfully reduce the error introduced in the perturbed landmarks and still produce accurate maps in the presence of perturbations.

\subsection{Benchmark test using the Faust Dataset} 

\begin{wrapfigure}[20]{r}{0.55\textwidth} 
\centering
\includegraphics[width=.95\linewidth]{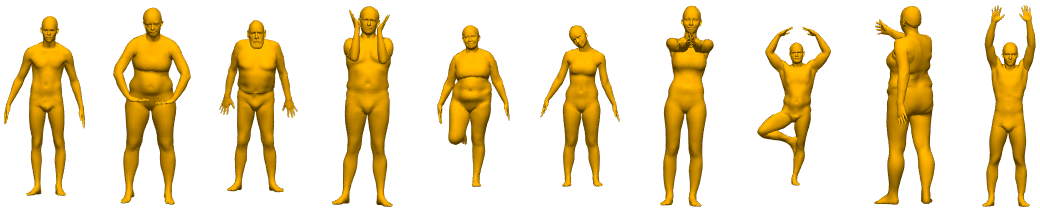}
\includegraphics[width=.9\linewidth]{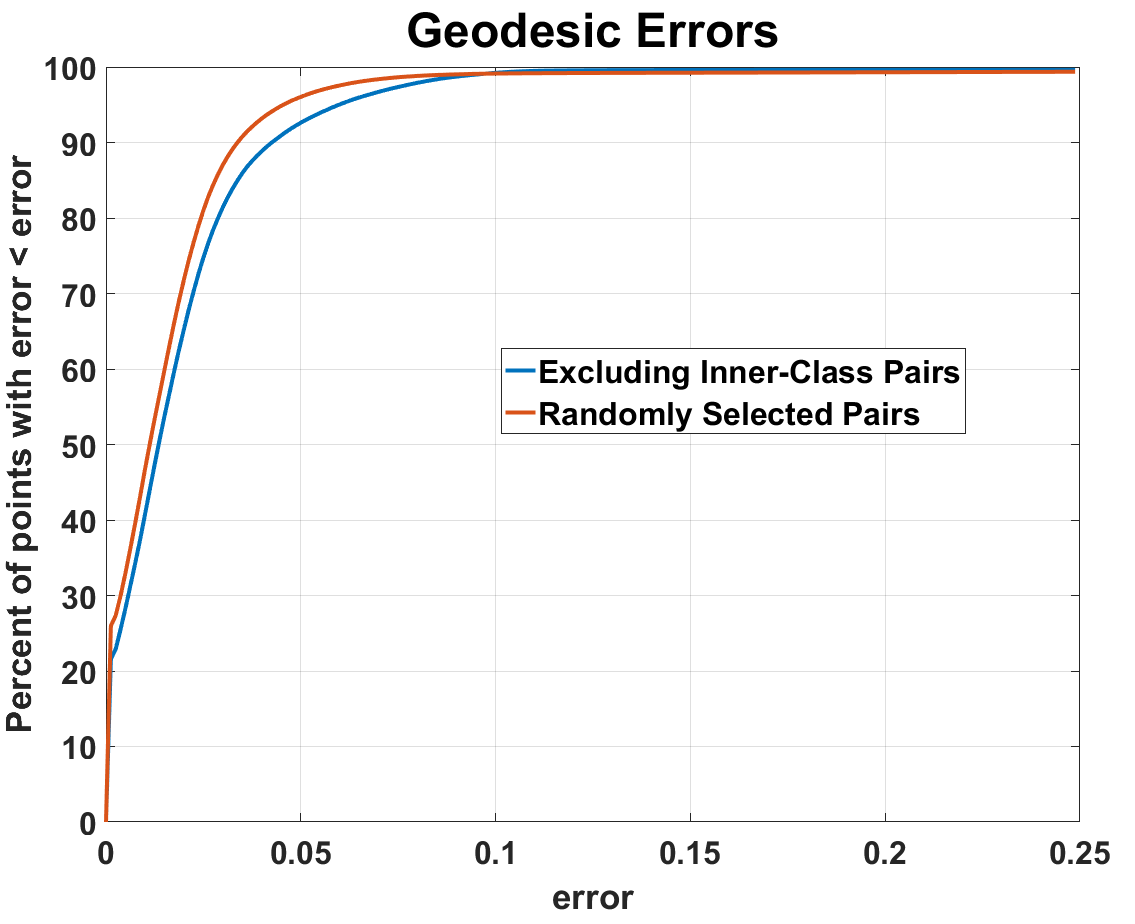}
\caption{Top: selected subjects from the Faust data set \cite{bogo2014faust}. Bottom: geodesic errors for randomly selected and least isomorphic pairs. }
\label{Faustex}
\end{wrapfigure}

In our next experiment, we test our algorithm on a larger dataset to demonstrate its effectiveness and robustness on a variety of shapes. The Faust dataset is a collection of 100 3D shapes composed of 10 real individuals in 10 distinct poses (Figure \ref{Faustex}) \cite{bogo2014faust}.
Instead of testing all 9900 possible correspondences be each of the pairs, we select two smaller subsets of shapes to formulate to smaller test sets. For the first test, we randomly choose 100 pairs of shapes and compute the correspondences. In the second test, we choose l0 scans and ensure that each individual and each pose is represented exactly once in the test set and compute all the 90 correspondences. This selection criteria ensures that no pairs are from the same the pose or individual. The bottom left graph in Figure \ref{Faustex} shows the average error of the mappings for each of these tests. We see that our algorithm again computes very accurate correspondences for both tests. Furthermore, we see that the results for the harder test set are very close to the results for the first test set. This indicates that our approach can effectively handle non-isometric matching problems with large deformations.

\section{Conclusions}
\label{sec:conclusion}
In this work, we have developed a variation method for computing correspondence between pairs of largely deformed non-isometric manifolds. Our approach considers conformal deformation of the manifolds and combines with traditional LB spectral theory. This method naturally connects metric deformations to the spectrum of the manifold and therefore allows us to register manifolds with large deformations. Our approach simultaneously aligns the bases of the manifolds and computes a conformal deformation without having to explicitly reconstruct the deformed manifolds. We have also proposed an efficient, locally convergent method to solve this model based on the PAM framework. Finally, we have conducted intensive numerical experiments to demonstrate the effectiveness and robustness of our methods.

\section*{Acknowledgement}
The research of S. Schonsheck and R. Lai is supported in part by NSF DMS--1522645 and an NSF Career Award DMS--1752934.


\end{document}